%
%
%
%
%

\input amstex
\documentstyle {amsppt}
\topmatter

\title
                ON THE REFLECTION PRINCIPLE IN ${\Bbb C}^n$
\endtitle
 \author
                        Alexander Sukhov
\endauthor
\abstract
We propose a reflection principle for holomorphic
     objects in ${\Bbb C}^n$ . Our construction generalizes
     the classical principle of H.Lewy, S.Pinchuk and S.Webster.

   Key-words: reflection principle, generic manifold, maximum
modulus set, holomorphic extension.
     A.M.S. Classification: 32D15, 32D99, 32F25, 32H99.
\endabstract
\endtopmatter
\document
\head
INTRODUCTION
\endhead
\par     The classical Schwarz  reflection principle asserts  that  a
conformal isomorphism between bounded domains ${D}\subset {\Bbb C }$
and  ${D}^{\prime}\subset{\Bbb C}$
with  real  analytic  boundaries  extends  holomorphically  to  a
neighborhood of  the  closure  $\overline D$.  In  several  variables,  this
phenomenon was first investigated by H.Lewy \cite{Le}  and  S.Pinchuk
\cite{Pi} for the  case  of  strictly  pseudoconvex  boundaries.  They
introduced a multidimensional reflection principle. Far  reaching
generalizations of their result were obtained by several  authors
( for instance, see \cite{BJT, DW, DF}, etc.).  Especially   we
would like to note the important paper  of  S.Webster  \cite{We1}.  He
discovered a new reflection  principle  that  differs  from  the
reflection principle of H.Lewy and S.Pinchuk. Using his technique
S.Webster  explored  the  extension  phenomenon  for  holomorphic
mappings between  domains  with  piecewise smooth  real  analytic
boundaries \cite{We2}.
\par     There is another way to develop the reflection principle  in
several  complex  variables.  Let  $ D $  and  $ {D}^{\prime}$    be
strictly pseudoconvex  domains  with  real  analytic   boundaries
in   ${\Bbb C}^n$ .
Suppose that $ f: D \to {D^{\prime}}$  is a  biholomorphic  mapping. One  can
assume that  $ f $  extends smoothly to the boundary $ bD$  of $ D $ by
Fefferman's theorem \cite{Fe}. Let $  A = {\Gamma}_f $  be the graph of $  f $   in
${\Bbb C}^{2n}$  .  Then $A$   is a complex n-dimensional
manifold with smooth boundary  $ bA \subset M$   and $ M = bD
\times {bD^{\prime}}$ is a generic real analytic
manifold. The extension theorem of H.Lewy  and  S.Pinchuk   means
that  $ A $  continues   to a complex manifold in a neidhborhood of
$M$.  Therefore,the following natural question arises. Let $\Omega$  be a
domain in ${\Bbb C}^N$   and let $ M$  be a generic real
analytic manifold in $\Omega$. Assume that $ A $ is a complex p-dimensional
manifold in ${\Omega} \backslash M$
with smooth boundary $bA \subset M $.  Under  what  conditions  does $A$
continue  analytically through $ M $?
\par H.Alexander  \cite{Al},   B.Shiffman   \cite{Sh}   and   E.Chirka
\cite{Ch1}
investigated  the  case  where $ A$  is  a   complex  1-dimensional
analytic set and  $ M $  is  a  totally  real manifold. In the general case an
affirmative  answer  depends  on the Levi form of
$M$. Also, a direction of the approach of $ A$   to  $M$  is essential.
In  this  connection  the  notion   of   Levi transversality  was
introduced \cite{PiCh,Su1,Su2,Su3}.   The condition  of  Levi
transversality  means   that   $ M $   has   a sufficiently  non -
degenerate Levi form  and  the  tangent   space   of  $A$   at   a
boundary point is  in the general position  with respect  to  the
Levi form of  $M$. Using the concept of Levi transversality one can
generalize  the  extension  theorem  of  H.Levi   and   S.Pinchuk
\cite{PCh,Su1,Su2,Su3}.
\par     However,this version of the reflection  principle  does  not
work if $ A $ is the graph of a holomorphic mapping between  wedges
with generic edges. Indeed, the condition $ bA \subset M $   presupposes
that ( locally )  $A$  is the graph of a holomorphic mapping over a
domain with smooth boundary. Hence, it seems natural to  consider
the case where the intersection  $\overline{A} \cap M $  has an  arbitrary
real dimension. The present paper is devoted to this problem.
\par     We shall consider a mapping   f   defined and holomorphic on
a wedge  $ W \subset {{\Bbb C}^p}$    with a smooth ( generic) edge  $E$.
  We  suppose that  $ f$   is smooth on $ W \cup E$. Let  $ A
\subset {{\Bbb C}^N}$    be the  graph  of
$f$   over $W$   and let  ${A}_{E}$    be the  graph  of  the  restriction
$f \vert E$   ( if $ E $ is a hypersurface, then  $ {A}_{E}  =  bA $ ).  We
suppose that  ${A}_{E}$  is contained in a generic real analytic
manifold $ M$.
\par     We introduce the notion of Levi transversality of $A$    and
$M$   at $ a \in{{A}_{E}} $  by analogy with \cite{Su2}.  This   condition   means
that  the Levi form of  $ M$   and the tangent space
${T}_{a}({A}_{E})$ are in
the  general position. The important example arises when   $A$   is
the  graph of a holomorphic mapping   $f$    between  wedges    $W$
and    $W^{\prime}$  with real analytic edges   $E$   and   $E^{\prime}$ ,  $
f(E)$ is
contained in $ E^{\prime}$
and   $M = E \times{E^{\prime}}$.   Our  first result  is  Theorem  1  that
establishes  the holomorphic extendability of   $f$   to  a neighborhood
of   $ E$   under the condition of Levi transversality. Moreover,
Theorem  1 asserts that in fact  $E$   is a real  analytic
manifold.  This theorem generalizes the above-mentioned results
on the  extension of holomorphic mappings and complex manifolds.
For instance,  the  graph of a  biholomorphism  $ f:  D \to {D^{\prime}}$
between  strictly pseudoconvex domains is Levi transverse to  $ M
=  bD\times{ bD^{\prime}}$;  a complex 1-dimensional manifold   $A$   with
smooth boundary  $ bA\subset M$    is Levi transverse to  $ M $  if $M$   is
totally  real,  etc. (see Propositions 2,3 and Corollaries 1-4 ).
\par     Our main tool is the generalized  Lewy - Pinchuk - Webster
reflection  principle  that was developed in \cite{PCh,Su1,Su2}.
We use  our  technique to push further an  interesting  recent
result  of  A.Nagel  and J.-P.Rosay \cite{NR} on maximum modulus sets
( subsets of the boundary of a domain $\Omega$   in  ${\Bbb C}^n$ , where a
holomorphic function  takes its maximum modulus ).  They proved
that if a maximum modulus set  $E$  is a smooth real  n-dimensional
submanifold in   $b{\Omega}$, and if   $b{\Omega}$   is strictly  pseudoconvex
and real analytic near  $E$, then  $E$   is real analytic (and
totally real). Our second result (Theorem 2) shows that both the
extension  theorem of H.Lewy and S.Pinchuk and the theorem  of
A.Nagel and J.-P.Rosay ( under the slight  additional  assumption
for   $E$   to be totally real  )  are  consequences of  a  general
reflection principle.
\par     All our considerations are purely local.
\par     This paper is organized  as  follows.  In  section 1  we  introduce
notations,  give  precise  definitions  and  statements  of   our
results. In sections 2-4 we develop the reflection principle  and  prove
Theorem 1. In section 5 we prove certain  corollaries  of  Theorem  1  (
Proposition 1, Corollaries 1, 4). In section 6 we prove Theorem 2.

\head
1. NOTATIONS, DEFINITIONS AND RESULTS
\endhead

\par     Let   $D$   be a domain in  ${\Bbb C}^p$    and  let $E \subset D$
be  a generic manifold of class  $ C^2$  and of real codimension
$ d$.
Then

  $$ E =  \lbrace {x \in D :  r_{j} (x) = 0,  j = 1,...,d.} \rbrace
\tag{1.1}$$
where the functions  $ r_{j} : D \to {\Bbb  R}$   are of class $C^2 (D)$
and
$\partial{r_1}\land...\land \partial{r_d} \neq 0$   on $D$.  Let $W(E,D)$   be
a wedge  with  the
edge   $E$:

          $$W(E,D) = \left \lbrace x \in D :  r_{j} (x) < 0 \rbrace \right.
\tag{1.2}$$

     Let  $G$  be a domain in  ${\Bbb  C}^n$ , and  $f: W(E,D) \to G$  be
a holomorphic mapping of class  $C^2$  on  $W(E,D) \cup E$.  Let us consider
the
graph  $A$   of $f$ :
   $$ A = {\Gamma}_{f}  =  \lbrace z = (x,y)  \in  {\Bbb   C}^{p}
\times {\Bbb C}^{n}  = {\Bbb  C}^N  : y = f(x), x \in W(E,D)
\rbrace , \tag{1.3}$$
that  is  a  complex  p-dimensional manifold  in  a  domain
${\Omega} = D \times G$.  We denote by   $ A_E$    the  graph  of  the
restriction
$f \vert E$:

   $$ A_E  = \left \lbrace z = (x,y): y = f(x), x \in E \rbrace \right.
\tag{1.4}$$

\par     Then  $ A_E $ is a real $( 2p - d )$-dimensional manifold  (  if
$d = 1$,   then  $ A_E  = bA = ( \overline{A} \cap {\Omega} ) \backslash A )$.

     Let $M$   be a generic real analytic  (closed)  manifold  of
real codimension $ m $  in a domain $ \Omega \subset {\Bbb  C}^N$ :

  $$ M = \lbrace z \in {\Omega} :  {\rho}_{j} (z,\overline{z}) = 0,
  j = 1,...,m \rbrace , \tag{1.5}$$
where the functions ${\rho}_{j} : {\Omega}\to{\Bbb  R}$  are
real  analytic
on
$\Omega$
and $  \partial{\rho}_{1}\land...\land \partial{\rho}_{m} \neq 0 $  on
$\Omega$.
We suppose that $ A_E\subset M$.  Our
main   question  :  under   what  conditions  does  $ A $  continue
analytically across  $M$ ?
\par     We  recall  certain  basic  definitions  of  the  theory  of
Cauchy-Riemann manifolds ( for instance, see \cite{Ch2}  ).  Let   $T_a M$
be the real tangent space of  $M$  at  $a\in M$. We  denote  by   $ T^c_a M $
the complex tangent space. We recall that  $ T^c_a M = T_a M \cap i(T_a
M)$.
     For  $ M$  of the form \thetag{1.5}  we have
$$
\aligned
& T_a M =  \lbrace t \in {{\Bbb  C}^N}  :\Re \sum^{N}_{{\nu} =
1}{{ \partial{\rho}_{j}} \over {\partial{z}_{\nu}}}(a) t_{\nu}  =
0,  j =1,...,m \rbrace, \\
&   T^c_a M =  \lbrace t \in {\Bbb  C}^N  : {\sum^{N}_{{\nu} =
1}}{{{\partial{\rho}_{j}}} \over {\partial{z}_{\nu}}}(a)t_{\nu}
= 0,  j =1,...,m\rbrace.\\
\endaligned
\tag{1.6}
$$
The complex dimension of  $ T^c_a M$   is equal to $ N - m$.  It is
called a CR dimension of $ M$. We denote the CR dimension  of   $ M $
by  $ CRdim M $ .
     For  $ u,v \in {\Bbb  C}^N$    and  $ a \in M$    we  denote  by
$H_{a}({\rho}_{j} ,u,v)$
the Levi form of $ {\rho}_{j}$  :

  $$   H_{a}({\rho}_{j} ,u,v) =
\sum^{N}_{{\nu},{\mu} = 1} {{{\partial}^2 {\rho}_{j}} \over
{{\partial z}_{\nu} {\partial \overline z}_{\mu}}}(a)u_{\nu}
{\overline v}_{\mu}.  \tag{1.7}$$

\par   Fix a hermitian scalar product $  < ,  >  $   in   ${\Bbb
C}^N$ .  We associate with each Levi  form  a hermitian
${\Bbb C}$-linear operator  $ L^{j}_{a} $
defined on $ T^{c}_{a} M$  by the condition
$ H_{a}({\rho}_{j},u,v) = < L^{j}_{a}(u),v> $  for any  $ u,v \in
T^{c}_{a}M$.  It is called the Levi operator of  $ {\rho}_{j}$ .

     Assume that  $ A_{E} \subset M$  and  $ m > d$.  We  shall  consider  the
case  when  the  restrictions  $ d{\rho}_{j} (a) \vert T_{a}A,  j
= m - d + 1,...,m$
are linearly independent ( perhaps , after a renumeration of  the
functions ${\rho}_{j}$   ).  Here  $ a \in A_{E} \subset  M$.   One  can
consider  this
requirement  as a condition of the  "partial"  transversality  of
the tangent spaces $ T_{a}A $ and  $T_{a}M$   (  taking  into  account  the
inclusion  $ A_{E} \subset M$ ).  Since  $ A_{E} \subset  M $  and
$ dim _{R} A_{E}  = 2p  -  d $,
the restrictions  $ d{\rho}_{j}(a) \vert T_{a}A, j = k,  m  -  d  +
1,...,m $   are
linearly dependent for any $  k = 1,...,m - d$.  Hence,  after  the
replacement of  ${\rho}_{k}$    by  ${\rho}_{k}  -
\sum^{m}_{j = m - d + 1} {\lambda}_{kj}{\rho}_{j}$  ,
one can assume that

  $$   d{\rho}_{k}(a) \vert T_{a}A = 0,  k = 1,...,m - d.
\tag{1.8}$$

     The following definition is a basic point of our approach.
\definition {Definition 1}
  Let  $ A $  and $  M $  be given as above. We say
that $  A $  and  $ M$   are Levi transverse at $  a  \in   A_{E}
\subset  M $  if the
following conditions hold:
\par     (i) $ m > d $;
\par     (ii) the restrictions  $ d{\rho}_{j}(a)|T_{a}A, j = m - d  +
1,...,m$  are
linearly independent ( perhaps,  after  a  renumeration  of  the
functions  ${\rho}_{j}$  );
\par     (iii) assume  that the functions ${\rho}_{j}$    are chosen  so  that
\thetag{1.8} holds. We suppose that

 $$      T^{c}_{a}(A_{E})   +    \sum^{m    -    d}_{j    =    1}
L^{j}_{a}(T^{c}_{a}(A_{E})) = T^{c}_{a}M. \tag{1.9}$$
\enddefinition
 {\it Remark 1}. We shall show that \thetag{1.9} does not depend  on  the
choice of a hermitian  scalar  product  that  defines  the  Levi
operators $ L^{j}_{a}$ .

   {\it  Remark 2}. Formally  one can treat  (i) as a consequence  of
(iii). Indeed, in the case $  m \leq d $ the condition (1.9)  has  the
form $ T^{c}_{a}(A_{E}) = T^{c}_{a}M$  or $ p  -  d  =  N  -  m$.
Since $  N -  p  >  0 $,  we
obtain (i). We prefer to impose (i) explicitely.

{\it Remark 3}. We emphasize that one requires (1.9) to be true at least
for a certain collection of the defining functions satisfying (ii) and (1.8)
(not necessarily for all such collections).

     Our first result is the following
\proclaim {Theorem 1}
 Let  $ {\Omega} \subset {\Bbb  C}^N  =
{\Bbb  C}^p \times {\Bbb  C}^n $    be a domain of the  form
${\Omega} = D \times G $,  where  $ D \subset {\Bbb  C}^p$    and
$  G \subset {\Bbb  C}^n$ .  Let  $  M \subset {\Omega}$  be  a
generic  real  analytic  manifold .  Suppose  that $  E \subset D
$ is  a  generic  manifold   of  class  $ C^2$    and  $ f: W(E,D)
\to G$   is a mapping holomorphic on the wedge  $W(E,D)$  and of
class  $ C^2$    on  $W(E,D) \cup E$.   Assume  that $ A$  ( = the
graph of  $ f $)   and  $ M $  are Levi transverse at a point $
\widetilde a = ( a,f(a)) \in A_E \subset M$.  Then   $f$ extends
holomorphically to a neighborhood of the point $ a \in  E $
and   $E$  is a real analytic manifold near $a$.
\endproclaim
\par     This theorem generalizes well-known results  connected  with
the reflection principle.  We  start  from  the  applications  of
Theorem 1 to the mapping problem.

     Let   ${\Omega}$   be a domain in ${\Bbb  C}^p$   and let   $  S
\in {\Omega}$   be a generic
real analytic manifold of the form

 $$ S = \lbrace x \in {\Omega} : r_{j}(x,{\overline x}) = 0,
j = 1,...,d \rbrace, \tag{1.10} $$
where the functions $  r_{j} : {\Omega} \to {\Bbb  R}$   are real
analytic on  $ {\Omega}$  and ${\partial  r}_{1}  \land  ...\land
{\partial r}_{d} \neq 0$   on  ${\Omega}$.  Let

 $$   W(S,{\Omega}) = \lbrace x \in {\Omega} :
 r_{j}(x,{\overline x}) < 0, j = 1,...,d \rbrace, \tag{1.11}$$
be a wedge in $ {\Omega}$   with the edge  $ S$. Similarly, let us consider
a domain  $ {\Omega}^{\prime} \subset {\Bbb  C}^{p^{\prime}}$    and a generic
real analytic manifold  $ S^{\prime} \subset {\Omega}^{\prime}$
of the form

 $$ S^{\prime} = \lbrace  x^{\prime} \in {\Omega}^{\prime} :
r^{\prime}_{j}(x,{\overline x^{\prime}}) = 0, j = 1,...,d^{\prime} \rbrace,
\tag{1.12}$$
where  the  functions  $  r^{\prime}_{j}:  {\Omega}^{\prime}  \to
{\Bbb  R}$   are real analytic  on  ${\Omega}^{\prime}$
and   ${\partial   r}^{\prime}_{1}   \land   ...\land    {\partial
r}^{\prime}_{d^{\prime}} \neq 0$   on  ${\Omega}^{\prime}$.   Fix  a
hermitian  scalar
product on  ${\Bbb  C}^{p^{\prime}}$ .  We denote  by $  L^{\prime
^{j}}_{a^{\prime}}$   the Levi operators  of  the
functions  $  r^{\prime}_{j}$   at $  a^{\prime} \in S^{\prime}$.
\par     The first consequence of Theorem 1 is the following
\proclaim {Proposition 1}
 Suppose that  $ {\Omega}$  is a domain
in ${\Bbb  C}^{p}$  , $ S $ is
a generic real  analytic  manifold  of  the  form  \thetag{1.10}  in   $
{\Omega}$,  ${\Omega}^{\prime}$
is a domain in ${\Bbb   C}^{p^{\prime}}$,  $  S^{\prime}  \subset
{\Omega}^{\prime}$  is a generic real analytic  manifold of the
form  \thetag{1.12}. Suppose  further  that $  f: W(S,{\Omega}) \to
{\Bbb  C}^{p^{\prime}}$     is  a  holomorphic mapping of class  $
C^2$   on  $ W(S,{\Omega}) \cup S $
and  $ f(S) \subset {S^{\prime}}$.  Assume that  for  certain   $a
\in  S $ the following condition holds:

 $$    \sum^{d^{\prime}}_{j = 1}
 L^{{\prime}^{j}}_{a^{\prime}}  (df_{a}(T^{c}_{a}S)) =
T^{c}_{a^{\prime}}(S^{\prime}), \tag{1.13}  $$
where  $ a^{\prime} = f(a)$   and  $ df_{a}$    is the tangent mapping. Then
$ f$ extends holomorphically to a neighborhood of $ a$.
\endproclaim
\par     In fact this statement is true for $ f \in C^1$ . This result  is
due to the author \cite{Su4,Su5,Su6}. A similar result (in another
form) was obtained by M.Derridj [Der]. We  shall  show  that  the
graph of $ f $ is Levi transverse  to  $ M = S \times  S^{\prime}
$.  Proposition 1
generalizes  the  classical   reflection   principle  of  H.Lewy,
S.Pinchuk  and S.Webster. To show this  fact  let  us  consider
the   Levi   form   $  Levi^{S^{\prime}}_{a^{\prime}}(u,v)$     of
$S^{\prime}$   at   $a^{\prime}$:

 $$    Levi^{S^{\prime}}_{a^{\prime}}(u,v) =
( H_{a^{\prime}}(r^{\prime}_{1},u,v),...,
H_{a^{\prime}}(r^{\prime}_{d^{\prime}} ,u,v) ). $$

     Recall, that manifold $ S^{\prime}$  is called Levi non-degenerate at
$ a^{\prime} \in  S^{\prime}$  if the equality
$  Levi^{S^{\prime}}_{a^{\prime}}(u,v)  =  0 $   for  any
$  v  \in  T^{c}_{a^{\prime}}(S^{\prime})$   implies  $ u = 0 $ ( see [We2] ).
\proclaim {Corollary 1}.
Suppose that $ {\Omega}$   is a domain in  ${\Bbb  C}^p$ ,  $ S
\subset {\Omega}$
is a generic real analytic manifold of the form \thetag{1.10},
${\Omega}^{\prime}$  is a
domain  in  ${\Bbb   C}^{p^{\prime}}$  ,   $   S^{\prime}   \subset
{\Omega}^{\prime}$  is a generic real  analytic  manifold  of
the form \thetag{1.12}. Suppose further that  $ f: W(S,{\Omega}) \to
{\Bbb  C}^{p^{\prime}}$    is  a
holomorphic  mapping  of class  $ C^2$    on
$  W(S,{\Omega}) \cup S $   and  $ f(S) \subset {S^{\prime}}$.
Assume that  the tangent mapping  $ df_{a} :T^{c}_{a} S \to
 T^{c}_{a^{\prime}}(S^{\prime})$  is
surjective and $ S^{\prime}$  is Levi non-degenerate at
$ a^{\prime} = f(a)$. Then   $ f $  extends holomorphically to  a
neighborhood of  $ a \in S$.
\endproclaim
\par     This result ( with  $ f $  of class  $ C^1$  ) is due to S.Webster
\cite{We2} ( see also \cite{TH} ).   We shall show in section 5 that  Corollary  1
is an immediate consequence of Proposition 1. We  would  like  to
emphasize that Proposition  1  is  a  considerably  more  general
statement. Indeed, the surjectivity of $ df_{a}$   implies
$ CRdim S^{\prime} \leq CR dim S $. From the other side \thetag{1.13} is valid
if     $ d^{\prime} CRdim S \geq CRdim S^{\prime}$. In particular, the
difference $ CRdim S^{\prime} - CRdim S $ can
be arbitrarily large.
\par     We point out that Proposition 1 and Corollary 1  deal with a
very special case of Theorem 1, because of in the  hypotheses  of
Theorem 1 $ M$  is not obliged to be the cartesian product.
\par     If both $ S$   and  $ S^{\prime}$  are  real  hypersurfaces  in
${\Bbb  C}^p$ ,  we
obtain from Corollary 1 the  classical  theorem  of   H.Lewy  and
S.Pinchuk.
\par     The next interesting special case of Theorem 1 arises if $  d
= 1$. Then  $ A_E  =  bA$,  the   condition   (i)   of   Definition  1
means  that   $  m \geq  2 $   and  (ii)  is  equivalent  to  the requirement
that $T_{0}A$  is not contained in  $T_{0}M$.  Thus,  we  obtain
the following
\proclaim {Proposition 2}.  Let $ M$  be a generic real  analytic
manifold
of real codimension $ \geq 2 $ in a  domain  $  {\Omega}  \subset
{\Bbb C}^N$ .  Suppose that  $ A $
is a complex  p-dimensional manifold  in   $  {\Omega} \backslash  M$   and
$(A,bA)$  is a $ C^2$ -manifold with boundary $ bA \subset M$. Assume that
$ A $ and $ M$  are Levi transverse at  $a \in  bA$. Then $ A $
continues    analytically
to  a  complex manifold in a neighborhood of $ a$.
\endproclaim
\par     This result was obtained by the author \cite{Su2}. In  fact  this
assertion is true for  $ (A,bA) \in C^1$   (see \cite{Su2}).  Proposition 2
generalizes  the  Lewy - Pinchuk extension  theorem as well.  We
obtain their statement  if  $ M $  is  the  cartesian  product  of
strictly  pseudoconvex  hypersurfaces  ${\Lambda}_{j} ,  j = 1,2$
and $ A$  is the graph of a mapping $  f: {\Lambda}_{1} \to
 {\Lambda}_{2}$   of class $ C^1$ ; $ f$  is holomorphic  on  a
domain $ D $  with  $ bD = {\Lambda}_{1}$    and  $ df $ is  non
-degenerate on  ${\Lambda}_{1}$ .
\proclaim {Corollary 1}   Let  $ M $  be  a  real  analytic  totally
real
N-dimensional manifold in a domain   $  {\Omega}  \subset  {\Bbb
C}^N$ .   Assume that  $ A \subset {\Omega} \backslash M$  is a complex
1-dimensional manifold and $ (A,bA) $ is a $ C^1$ -manifold  with
boundary   $ bA \subset  M$.  Then   $ A $    continues
analytically to a complex manifold in a neighborhood of $ bA$.
\endproclaim
\par     For complex 1-dimensional analytic sets similar results were
obtained
by H. Alexander \cite{Al}, B.Shiffman \cite{Sh} and E.Chirka \cite{Ch1}.
One  can  consider  Corollary  2  as  a  generalization  of   the
classical Schwarz  reflection principle. Indeed, if $  f:  D  \to
D^{\prime}$  is  a  conformal  isomorphism  between  domains   in
${\Bbb C}$   then the graph of  $ f $  is a  complex
1-dimensional
manifold whose boundary is
contained  in  the  totally  real  manifold  $   M  =  bD  \times
bD^{\prime} \subset {\Bbb  C}^2 $ .
\proclaim {Corollary 3}. Let  $ M $  be a generic real  analytic
manifold
in a domain  ${\Omega} \subset  {\Bbb   C}^N$     and   let   $  A
\subset  {\Omega} \backslash  M $  be  a  complex
p-dimensional  manifold;  we   suppose   that   $ (A,bA)$    is   a
$C^1$ -manifold with boundary $ bA \subset  M$.   Assume  that  $
CRdim M = p - 1$.
Then  $ A$   continues  to a complex  manifold in a neighborhood of
$bA$.
\endproclaim
\par     E.Chirka \cite{Ch1} obtained a similar result for analytic sets,
see also \cite{Fo}.
\par     Theorem 1 also implies the next
\proclaim {Corollary 4}.  Let   ${\Omega} \subset {\Bbb C}^N  =
{\Bbb  C}^p \times  {\Bbb  C}^n$    be  a  domain  of  the
form  $ {\Omega} = D \times G$,  where  $ D \subset {\Bbb  C}^p$
and   $G \subset {\Bbb  C}^n$ . Let  $ M \subset {\Omega}$  be
a  real analytic totally   real   N-dimensional  manifold,   $  E
\subset D$
be a totally  real   p-dimensional  $ C^1$ -manifold .  Suppose that
$f: W(E,D) \to  G $  is  a  holomorphic mapping  of class  $ C^1$    on
$W(E,D) \cup E$.  Assume that $ A_E  \subset M $( here $ A_E$   is
the graph of the
restriction  $f \vert E$ ).  Then   $f$   extends  holomorphically  to  a
neighborhood of  $ E$   and $ E$   is a real analytic manifold.
\endproclaim
\par     We shall use Corollary 4 to prove our second main result.
\proclaim {Theorem 2}.  Let  ${\Omega}$   be a domain in  $ {\Bbb
C}^n$    and let  $ E \subset  b{\Omega}$
be a totally real  n-dimensional $ C^1$ -manifold. Assume that  $
b{\Omega}$ is
strictly pseudoconvex  and  real  analytic  near  $  E$.  Let   $
E^{\prime}$   be  a
totally real $  n^{\prime}$-dimensional  real  analytic  manifold  that  is
contained in the boundary of a domain  $ {\Omega}^{\prime} \subset
{\Bbb  C}^{n^{\prime}}$.  .  Suppose  that
for  a  neighborhood    $V  \supset   E^{\prime}$    in   ${\Bbb
C}^{n^{\prime}}$    one has $ {\Omega}^{\prime} \cap V =
\lbrace x^{\prime} \in V: \rho(x^{\prime},{\overline x}^{\prime}) <
0 \rbrace $,  where   ${\rho}$   is a real analytic plurisubharmonic
function on  $V$,  $d{\rho} \neq 0$. Assume that $  U $  is a
neighborhood  of  $E$  in ${\Bbb  C}^n$   and
$ f: {\Omega} \cap U \to {\Omega}^{\prime} \cap V$   is a
holomorphic mapping of class   $C^2$    on $  {\overline{\Omega}}
\cap U$  such that $  f(E) \subset E^{\prime}$. Then  $ f $  extends
holomorphically to a neighborhood  of  $ E$   and  $ E$   is  a  real
analytic manifold.
\endproclaim
\par     We would like to  emphasize  that  we  do  not  require  the
inclusion   $ f(b{\Omega}) \subset  b{\Omega}^{\prime}$.  Theorem  2  implies
the   following
well-known assertion that is a version of the extension  theorem
of H.Lewy and S.Pinchuk.
\proclaim {Corollary 5}   Let    ${\Omega}$       and
${\Omega}^{\prime}$  be  bounded  strictly pseudoconvex domains
with  real  analytic  boundaries  and  let  $  f:  {\Omega}   \to
{\Omega}^{\prime}$  be  a biholomorphic mapping. Then  $ f$
extends  holomorphically   to   a  neighborhood  of  the  closure
$\overline {\Omega}$.
\endproclaim
\par     Indeed, Fefferman's theorem \cite{Fe} implies  that  $ f \in
  C^{\infty}(\overline{\Omega})$
and  $ f: \overline{\Omega} \to \overline{\Omega}^{\prime}$    is  a
diffeomorphism.  Let us  consider  a
point  $a \in b{\Omega}$  and a totally real  n-dimensional real
analytic
manifold  $ E^{\prime}  \subset  b{\Omega}^{\prime}$  through
$f(a)$. . Then $  f^{-1}(E^{\prime}) = E \owns a $
is   a   smooth   totally   real    n-dimensional   manifold   in
$b{\Omega}$  and one can apply Theorem 2.
\par     Another  application  of  Theorem  2  considers      maximum
modulus sets. Let   ${\Omega}$    be  a  strictly  pseudoconvex  (bounded)
domain in $ {\Bbb  C}^n$    with $ C^2$ -boundary. A subset
$  E \subset b{\Omega}$   is called a maximum modulus set  if  for
every  $a \in E$ there  exists  a  neighborhood  $U$ of $a$ and a
function  $ f$   defined and continuous on
$\overline{\Omega} \cap  U$,  holomorphic  on
${\Omega} \cap U$  so that $ \vert f \vert < 1$  on
${\Omega} \cap U $ and  $ \vert f \vert \equiv  1$  on  $ E  \cap
U$.
\par     In the case the function  $ f$   in the definition of  maximum
modulus set can be chosen of class  $ C^k$    on
${\overline {\Omega}} \cap  U$,  $  E $  is
called a  $C^k$ -maximum modulus set ( see \cite{NR} ).
\par     Theorem 2 immediately implies the following
\proclaim {Corollary 6} Let  ${\Omega} \subset {\Bbb  C}^n$   be a
domain and let  $E \subset b{\Omega}$  be  a
totally real n-dimensional manifold  of   class   $  C^1$  .   If
$b{\Omega} $  is
strictly pseudoconvex and real analytic near $ E$,  and  if $   E
$ is a  $C^2$ -maximum modulus set, then $ E$  is real analytic.
\endproclaim
\par     This result was obtained by A.Nagel and J.-P.Rosay \cite{NR}.  In
fact they proved it without the assumption that $  E $  is  totally
real. This requirement is not essential  for  our  proof  and  we
claim it for the convenience.
\par     There is the firm confidence that Theorem 1 is true if both
$f$  and $ E$  are of class $ C^1$ . To avoid  technical   complications
we consider here only the  $C^2$   case.

\head
2. BASIC PROPERTIES OF THE LEVI TRANSVERSALITY.
\endhead

     The concept of Levi transversality  is  obviously  invariant
under    biholomorphic    transformations.    Hence,     without
loss of generality one can assume that $  A$   and $  M$    are  Levi
transverse at  $ 0 \in A_E \subset  M$.
\par     The  condition   of   Levi   transversality   imposes   the
restrictions on the dimensions of   $ A$   and   $ M$.   Indeed,  the
inclusion   $ A_E \subset  M$   implies

    $$ dim T^c_{0}(A_E) = p - d \leq  dim T^c_{0}M = N - m $$.

     From the other side we have

  $$    dim \lbrack  T^c_{0}(A_E) +
    \sum^{m - d}_{j = 1} L^{j}_{0}(T^{c}_{0}(A_E)) \rbrack
\leq  (m - d + 1)(p - d) $$,

and \thetag{1.9} implies

              $$     (m - d + 1)(p - d) \geq N - m. $$

\proclaim {Lemma  2.1} There    exists    a
non-degenerate
${\Bbb C}$-linear change     of   coordinates  in  $  {\Bbb   C}^N$
such   that   in  the  new coordinates ( perhaps,  after  a
replacement  of  the  defining
functions by their linear combinations ) we have

  $$ \rho_{j}(z,{\overline z}) =
 z_{N - m + j} + {\overline z}_{N - m + j}  + o(\vert z \vert),
j = 1,...,m, $$

   $$  T_{0}A = \lbrace t \in {\Bbb  C}^N  :
t_{p - d + 1} =...= t_{N - d}    = 0 \rbrace, \tag{2.1}$$

$$      T^{c}_{0}(A_E) = \lbrace t \in {\Bbb  C}^N
  : t_{p - d + 1} =...= t_{N} = 0 \rbrace.$$

     The condition \thetag{1.9} holds for the new defining functions.
\endproclaim
\par    {\it Proof}. Since  $ A_E  \subset M$   and
$ d{\rho}_{m - d + 1} \land ...\land  d{\rho}_{m}  \vert   T_{0}A
\neq 0 $,  we have   $ T_{0}(A_E) = T_{0}A \cap T_{0}M$.  Also, we
have   $ d{\rho}_{j}(0) \vert  T_{0}A = 0, j = 1,...,m  -  d $.
Therefore,  taking  into  account  the  equality
$ dim T^{c}_{0}(A) =  p  -  d$,   we  conclude  that   $\partial
{\rho}_{m - d + 1} \land ...  \land  \partial{\rho}_{m}  \vert
T_{0}{A} \neq 0$. After  a  non-degenerate  ${\Bbb  C}$-linear   change  of
coordinates  we obtain

  $$      T^{c}_{0}M = (e_{1} ,...,e_{N - m})_{\Bbb  C} ,
T^{c}_{0}(A_E) = (e_{1},...,e_{p - d})_{\Bbb  C}  $$
where  $ e_{j} , j =  1,...,N$   is   the   standard   basis   in
${\Bbb  C}^N$    and  $(   )_{\Bbb  C}$   is the ${\Bbb C}$-linear hull.
Since  $ T^{c}_{0}(A_E) = T^{c}_{0}M \cap T_{0}A$,   we
have  $T_{0}A = T^{c}_{0}(A_E) \oplus  (v_{1},...,v_{d}  )_{\Bbb
C}$   and the vectors  $ e_{1},...,e_{N - m}, v_{1} ,...,v_{d} $    are
linearly independent. Let us consider  $v_{1} ,...,v_{d}$
as the vectors   $e^{\prime}_{N - d + 1},...,e^{\prime}_{N}$   of a
new basis in ${\Bbb  C}^N$   (  we  do
not change other basis vectors ). Then ( we omit the primes )

 $$    T_{0}A = (e_{1},...,e_{p - d})_{\Bbb  C}  \oplus L,
  L = (e_{N - d + 1},...,e_{N})_{\Bbb  C} . $$

     One  can  represent  $  T_{0}M $   in  the  form   $ T_{0}M   =
T^{c}_{0}M \oplus V$, where  $V$  is the othogonal complement  of
$ T^{c}_{0}M $ in $ T_{0}M$   with respect to
the standart real scalar product on   $ {\Bbb  C}^N  \cong {\Bbb
R}^{2N}$  .  Then  $ V \subset {\Bbb  C}^m  =
(e_{N - m + 1},...,e_{N})_{\Bbb  C}$.    Moreover,   $ dim_{\Bbb
R} V = m $,  $  ( V )_{\Bbb  C}  = {\Bbb  C}^m $,
$dim_{\Bbb  R} V \cap L = d$  and  $( V \cap L )_{\Bbb  C}  = L$.  One can
represent ${\Bbb  C}^m$   in the form $ {\Bbb  C}^m  = L + P$,  where
$ P$  is a complex $(m - d)$-dimensional linear space,  $V =  (V
\cap L) \oplus (V \cap P)$, $dim_{\Bbb  R} V \cap P = m - d $
and   $ (V \cap P)_{\Bbb  C}  = P $.  After a non-degenerate
${\Bbb C}$-linear  change
of coordinates in  ${\Bbb  C}^m $  we obtain  $ L = {\Bbb  C}^d$ ,
$ P = {\Bbb  C}^{m - d} $  , $  V \cap L  = i{\Bbb  R}^d  $,
$ V \cap P = i{\Bbb  R}^{m - d}$
and   $V = i{\Bbb  R}^d  + i{\Bbb  R}^{m - d}  =  i{\Bbb  R}^m$ .
Thus,  we have  $ L = (e_{N - d + 1},...,e_{N})_{\Bbb  C}$    and
$  V = i{\Bbb  R}^m \subset (e_{N - m + 1},...,e_{N})_{\Bbb  C}$ .
Then  $ T_{0}M = \lbrace t \in {\Bbb  C}^N  : t_{j}  +
{\overline t_{j}}  = 0, j = N - m + 1,...,N \rbrace $. After
a replacement of the functions $ {\rho}_{j}$     by  their (real) linear
combinations ,   we obtain

  $$   {\rho}^{\prime}_{j} = z_{N - m + j}  + {\overline z_{N - m
+ j}}  + o(\vert z \vert),  j = 1,...,m.$$

     Since  $ d{\rho}_{m - d + 1} \land ... \land d{\rho}_{m}  \vert
T_{0}A \neq 0  $,   we  have  ${\rho}^{\prime}_{j}  =   \sum^{m  -
d}_{{\nu} = 1} {\lambda}_{j{\nu}}{\rho}_{\nu}$
for $ j = 1,...,m - d $  and  $ det({\lambda}_{j{\nu}}) \neq 0$. Hence,
for any complex
linear  space  $  S  \subset T^{c}_{0}M$    we obviously have
$ \sum^{m - d}_{j = 1} L^{{\prime}^j}_{0}(S) = \sum^{m - d}_{j  =
1}   L^{j}_{0}(S)$   ,    where    $L^{{\prime}^j}_{0}   $    and
$L^{j}_{0}$   are the Levi  operators  of  ${\rho}_{j}^{\prime} $  and
${\rho}_{j}$   respectively,{ \it Q.E.D.}
\par    It is convenient to represent $  A$   as a graph over $ T_{0}A
\cong {\Bbb  C}^p$  in the coordinates \thetag{2.1}. We set
$ z^{\prime} = (z_{1},...,z_{p - d})$, $ z^{\prime \prime} =
(z_{p - d + 1},...,z_{N - d})$, $ z^{\prime \prime \prime} =
(z_{N - d + 1},...,z_{N})$.
\proclaim {Lemma 2.2} Assume that the  coordinates  \thetag{2.1}
are chosen. Then ( in a neighborhood of the origin )  $ A $  can
be represented
in  the  form  $  z^{\prime  \prime}   =   g(z^{\prime},z^{\prime
\prime \prime})$.  Here  $ g $   is  a  mapping
holomorphic on a wedge  $ W(\widetilde E, \widetilde D) \subset
{\Bbb  C}^p  = T_{0}A $  with a  $ C^2$ -edge  $\widetilde E$
and $  g $   is  of  class   $  C^2$     on   $  W(\widetilde  E,
\widetilde D) \cup {\widetilde E}$. The  graph  of  the
restriction $ g \vert {\widetilde E}$  coincides with  $ A_E$ .
\endproclaim
\par     This is a  simple  consequence  of  the  implicit  function
theorem. We leave details to the reader.
\par     In what follows we omit the tilde and assume that

$$     E = \lbrace (z^{\prime},z^{\prime \prime \prime})
\in D : r_{j}(z^{\prime},z^{\prime \prime \prime}) = 0,
 j = 1,...,d \rbrace, $$

$$    W(E,D) = \lbrace (z^{\prime},z^{\prime \prime \prime})
\in D :r_{j}(z^{\prime},z^{\prime \prime \prime})
< 0, j = 1,...,d \rbrace. $$
\proclaim { Lemma 2.3}
 One has $ H_{0}({\rho}_{j} ,u,u) = 0 $ for  $ u \in
T^{0}_{0}(A_E)$   and $j = 1,...,m - d$.
\endproclaim
\par   {\it  Proof}. Since  $ A_E  \subset M$,   we have
${\rho}^{j}     =    {\rho}_{j}(z^{\prime},g(z^{\prime},z^{\prime
\prime        \prime}),z^{\prime        \prime         \prime},
{\overline z^{\prime}},{\overline g(z^{\prime},z^{\prime \prime \prime}}),
{\overline  z^{\prime  \prime
\prime}}) = 0$,  for $ j = 1,...,m$   and   $(z^{\prime},z^{\prime
\prime \prime}) \in E$.  By the division lemma we  obtain
${\rho}^{j} = \sum {\lambda}_{j{\nu}}r_{\nu}(z^{\prime},z^{\prime
\prime \prime})$,   where the functions  ${\lambda}_{j{\nu}}$     are
of class  $ C^1$    in a neighborhood
of the origin in  ${\Bbb  C}^p$ .   Since
$ d{\rho}_{j}(0) \vert T_{0}A = 0,  j = 1,...,m - d $,
we conclude that $  {\lambda}_{j{\nu}}(0) = 0 $  for
$  j = 1,...,m - d,  {\nu} = 1,...,d$. Hence,
${\rho}^{j}  = o({\vert z^{\prime} \vert}^2  + {\vert z^{\prime \prime
\prime} \vert}^2)$   for $ (z^{\prime},z^{\prime \prime \prime}) \in
 T^{c}_{0}E$. This implies our statement.{\it Q.E.D}.
\par     Let us show that \thetag{1.9} does not depend on the choice of  the
hermitian scalar product which defines the Levi operators.  Let
$< , >^{\prime}$  be another hermitian scalar  product  and  let
the  Levi operators  $L_{0}^{\prime^{j}}$   be defined by the condition
$<L^{\prime}_{j}(u),v>^{\prime} = H_{0}({\rho}_{j},u,v) $  for
$u,v \in T^{c}_{0}M$.  Lemma 2.3 implies

 $$     H_{0}({\rho}_{j},u,v) = < L^{j}_{0}(u),v > =
< L^{\prime^{j}}_{0}(u),v>^{\prime} = 0 $$
for $  j = 1,...,m - d $  and  $ u,v \in T^{c}_{0}(A_E )$.  Also,
it is easy  to
show that $ L^{\prime^{j}}_{0}   =  R L^{j}_{0}$ ,     where $  R
$   is  a  non-degenerate ${\Bbb C}$-linear operator on $ T^{c}_{0}
M$.  Hence, the spaces $  S = \sum^{m - d}_{j = 1}
L^{j}_{0}(T^{c}_{0}A_E) $
and      $    S^{\prime}    =    \sum^{m    -    d}_{j    =    1}
L^{\prime^{j}}_{0}(T^{c}_{0}A_E)$
  have the same dimension. But $ S $ (resp. $S^{\prime}$)  is contained in
the orthogonal complement of $ T^{c}_{0}(A_E)$   with
respect to $ < , >$  ( resp.$ < , >^{\prime}$ ). Thus, if \thetag{1.9} holds for
the operators $ L^{j}_{0}$   then \thetag{1.9} also holds for
$L^{\prime^{j}}_{0}$ .
\par     Now we are going to investigate the condition \thetag{1.9}  in  the
coordinates \thetag{2.1}. Without loss of generality one can assume that
$< z,   w   >  =  \sum  z_{j}{\overline  w_{j}}$  in  the  coordinates
\thetag{2.1}
and this  hermitian scalar product  defines  the  Levi  operators
$L^{j}_{0}$  .
\par     For  positive  integers  $ j, {\nu}_{1} ,...,{\nu}_{k}$     we
consider   the matrices of the following form

 $$ H(j,{\nu}_{1},...,{\nu}_{k}) =
{ \left ({{\partial^{2}{\rho_{j}}} \over {{\partial
z_{\nu}}{\partial{\overline  z_{\mu}}}}}(0)  \right )  } ^{{\nu}   =
{\nu}_{1},...,{\nu}_{k}}_{{\mu} = p - d + 1,...,N - m},
\tag{2.2}$$

( there are $ k$  rows and $ (N - m - p + d)$  columns in \thetag{2.2}).
\proclaim { Lemma 2.4} Assume that the  coordinates  \thetag{2.1}  are  chosen.
Then the condition of Levi transversality is  equivalent  to  the
following:  there exist  collections $ (j(n),
{\nu}_{1}(n),...,{\nu}_{k(n)}(n)) $, $n = 1,...,s$   of positive
integers so that  $ 1 \leq j(1) <...<  j(s) \leq m - d $, $  1 \leq
{\nu}_{1}(n)  <...< {\nu}_{k(n)}(n) \leq p  -  d$    for   $  n  =
1,...,s$, $\sum^{s}_{n = 1} k(n) = N - m - p + d$   and

           $$                 det H \neq 0, \tag{2.3}$$
where  $ H $  is  a   $ (N - m - p + d) \times (N - m - p +
d)$-matrix  of
the form

$$    H =  \pmatrix
H(j(1),{\nu}_{1}(1),...,{\nu}_{k(1)}(1)) \cr
\hdots \cr
H(j(s),{\nu}_{1}(s),...,{\nu}_{k(s)}(s))
\endpmatrix . \tag{2.4}$$
\endproclaim
   {\it  Proof}.  Let $ e_{j} , j = 1,...,N$  be the standard basis  in
${\Bbb  C}^N$ . Since the spaces $  L^{j}_{0}(T^{c}_{0}A_E)$   are
contained  in  the  orthogonal
complement of  $  T^{c}_{0}A_E$  (  see  Lemma  2.3  ),  the condition
\thetag{1.9}  is
equivalent to

$$  dim \lbrack \sum^{m - d}_{j = 1}
L^{j}_{0}(T^{c}_{0}A_E) \rbrack   = N - m - p + d$$.

     Hence, there are  $( N - m - p +  d )$  ${\Bbb C}$-linearly independent
vectors among the vectors $  L^{j}_{0}(e_{k}), j = 1,...,m - d, k =  1,...,p
- d$.
\par     Let  $ H_{j}$    be the matrix of the restriction of
$  H_{0}({\rho}_{j},u,u)$ on  $T^{c}_{0}M$   ( with respect to the basis
$ e_{\nu} , {\nu} = 1,...,N - m )$.
Each Levi operator has the matrix ( with respect to this basis )
that is transpose  of $ H_{j}$ . But the coordinates of the vector
$ L^{j}_{0}(e_{k})$ form the $ {\nu}$-th column of the matrix of
$ L^{j}_{0}$ .  Let us consider  the
system of vectors that is the union over $ j = 1,...,m  -  d $   of
the first $ ( p - d )$  rows of $ H_{j}$ .  It was just shown that  there
are  $ ( N - m - p + d)$  ${\Bbb  C}$-linearly independent vectors  in  this
system.  Lemma 2.3  implies that

$$         {{\partial^{2}{\rho}_{j}}         \over          {{\partial
z_{\nu}}{\partial{\overline z_{\mu}}}}}(0) = 0 , {\nu},{\mu} = 1,...,p
- d, j = 1,...,m - d. $$

     We  consider  the  matrix  $  C$    which  is  formed  by  the
above-mentioned  $( N - m - p + d )$  independent rows.   Then  $
C$  has the form  $ C = ( O \vert H )$,  where $  H$   is the  matrix
\thetag{2.4}.
Hence, the linear independence of the rows of $ C $  is  equivalent
to (2.3), {\it Q.E.D}.
\head
3. COMPLEXIFICATION.
\endhead

    Assume the coordinates \thetag{2.1} are chosen. According  to  Lemma
2.2 we have in a neighborhood  $ U \owns 0$ such that

  $$            A \cap U = \lbrace z \in U : z^{\prime \prime} =
g(z^{\prime},z^{\prime \prime \prime}) \rbrace, $$
where  $ g $ is a mapping holomorphic on $ W(E,D)$  and of class
$  C^2$   on  $ W(E,D) \cup E$. Moreover,$ g(0) = 0$  and  $dg(0)
= 0$.
     We set

 $${\rho}^{j}(z^{\prime},z^{\prime \prime \prime}) =
( {\rho}_{j} \vert A ) = {\rho}_{j}(z^{\prime},
g(z^{\prime},z^{\prime \prime \prime}),z^{\prime \prime \prime},
{\overline z}, {\overline g(z^{\prime},z^{\prime \prime \prime})},
{\overline z^{\prime \prime \prime}}).$$

     Since  $ A_E  \subset M$,  one has  $ {\rho}^{j} \vert  E  =
0$    for $  j  =  1,...,m$.
By   $ d{\rho}^{m - d + 1} \land... \land d{\rho}^{m}  \neq 0 $,   we get

$$   E = \lbrace (z^{\prime},z^{\prime \prime \prime})
\in {\Bbb  C}^p  : {\rho}^{j}  = 0,
  j = m - d + 1,...,m \rbrace, \tag{3.1}$$
in a neighborhood of the origin.
\par     Let us consider the vector fields  $ T_{q} ,  q  =   1,...,p  -  d$
( the sections of the complex tangent bundle  $ T^{c}E$)  of the form

 $$   T_{q}  = {\partial{} \over {\partial z_{q}}} -
 \sum^{d}_{j = 1} a_{jq}(z^{\prime},z^{\prime \prime \prime},
{\overline z^{\prime}},{\overline z^{\prime \prime \prime}})
{\partial{} \over {\partial z_{N - d +j}}}, \tag{3.2}$$
where

 $$   a_{jq} = \sum^{d}_{s = 1} b_{js} {{\partial {\rho}^{s + m -
d}} \over {\partial z_{q}}}. \tag{3.3} $$

     Here  $ ( b_{js})$   is the inverse of the Jacobian matrix

$$ S =  {\left ( {{\partial{\rho}^{k}} \over {\partial  z_{l}}}
\right )} ^{k = m - d +1,...,m}_{l = N - d + 1,...,N}.
\tag{3.4}$$

     It is easy to see that the vector fields  $ T_{q} , q = 1,...,p -
d$   form a basis of the bundle $  T^{c}E $  over a neighborhood of the
origin in  $ E$.
\proclaim {  Lemma 3.1}         For  $  (z^{\prime},z^{\prime  \prime
\prime}) \in E$  we have

$$    T_{\nu}{\rho}^{j}  = 0 ,  {\nu} = 1,...,p - d, j =  1,...,m
- d $$
in a neighborhood of the origin.
\endproclaim
\par   {\it  Proof}. Since   ${\rho}^{j} \vert E  =  0$,   there   are
$C^1$  -functions    ${\lambda}_{jk}(z^{\prime},z^{\prime  \prime
\prime})$   such that

$${\rho}^{j}  = \sum^{m}_{k = m - d + 1} {\lambda}_{jk}{\rho}^{k} ,
j = 1,...,m - d, $$
where  $ (z^{\prime},z^{\prime \prime \prime})$   belongs to a
neighborhood of  the  origin  in  ${\bold C}^p$ .  Hence, for
$(z^{\prime},z^{\prime \prime \prime}) \in E$
we have

$$ {{\partial{\rho}^{j}} \over {\partial z_{q}}}  =
\sum^{m}_{k = m - d + 1} {\lambda}_{jk} {\partial {\rho}^{k}
\over {\partial z_{q}}}, \tag{3.5}$$
for  $ j = 1,..., m - d, q = 1,...,p - d, N - d + 1,...,N$.
\par      Fix an arbitrary  $ j$. We assume that $ q = N - d + 1,... ,
N$  in  \thetag{3.5}  and consider \thetag{3.5} as a system  of linear equations
in ${\lambda}_{jk},  k = m - d + 1,...,m$.
     Then the matrix of this system is transpose of $  S$  of the
form \thetag{3.4}. Let us consider the rows

$$  U_{q}  = \left (
 {{\partial{\rho}^{m - d + 1}} \over
{\partial z_{q}}} \hdots {{\partial{\rho}^{m}} \over {\partial z_{q}}}
\right ), \tag{3.6}$$

and the columns

$$ {\Lambda}_{j} = \pmatrix
{\lambda}_{j m - d + 1} \cr \vdots \cr {\lambda}_{jm}
\endpmatrix,
V_{j} = \pmatrix
{{\partial {\rho}^{j}} /     {\partial z_{N - d + 1}}} \cr
 \vdots \cr {{\partial {\rho}^{j}} /     {\partial z_{N}}}
\endpmatrix. \tag{3.7} $$

      Since the matrix \thetag{3.4} is non-degenerate in the origin, one
has  $ {\Lambda}_{j}   =  (S)^{t^{-1}}  V_{j}$  ,    where   $  (
)^t$  denotes  the  transposed
matrix. Then \thetag{3.5} implies

$$  {{\partial {\rho}^{j}} \over {\partial z_{q}}}
= U_{q}{\Lambda}_{j} = U_{q}(S)^{t^{-1}}V_{j},  $$
for $ j = 1,...,m - d,  q = 1,...,p - d$.
     Hence,

$$  {{\partial{\rho}^{j}} \over {\partial z_{q}}} -
(V_{j})^{t}(S)^{-1}(U_{q})^{t} = 0 $$

     In view of \thetag{3.2}, \thetag{3.3},
\thetag{3.4},\thetag{3.6},\thetag{3.7} the  last  equality
is exactly equivalent to the assertion of Lemma 3.1, {\it Q.E.D}.
\par     Let   us    consider    the    functions    $     \widetilde
{\rho}_{j}(z,{\zeta})$    that  we  obtain
after the substitution  of ${\zeta} \in {\Bbb  C}^N$    instead  of
${\overline z}$   in  the
expansion of  ${{\rho}_{j}}(z,{\overline z})$    at  the
origin.  Then  the  functions
$\widetilde {{\rho}_{j}}(z,{\zeta})$   are holomorphic on a neighborhood of
the origin in  ${\Bbb  C}^{2N}$
and  $\widetilde {{\rho}_{j}}(z,{\zeta}) \vert
\lbrace  \zeta = \overline z \rbrace =  {\rho}_{j}(z,{\overline
z})$. We set

$$ M^{c}  = \lbrace (z,{\zeta}) \in {\Bbb  C}^{2N}   :
{\widetilde {{\rho}_{j}}}(z,{\zeta}) = 0, j = 1,...,m \rbrace.
\tag{3.8}$$

     Then $ M^{c}$    is a complex $ (2N - m)$-dimensional manifold in a
neighborhood  of  the  origin   in ${\Bbb  C}^{2N}$  .   It   is   called   a
complexification of $ M$. Also,  we introduce a real manifold $
{\widehat A_E}  = \lbrace (z,{\overline z}) \in {\Bbb  C}^{2N}
: z \in A_E  \rbrace$.   Then  ${\widehat A_E}  \subset M^{c} \cap
\lbrace {\zeta} = {\overline z} \rbrace$    and
$ dim_{\Bbb  R} {\widehat A_E}  = 2p - d$.
\par     Set $ \widetilde {{\rho}^{j}}(z^{\prime},z^{\prime  \prime
\prime},{\zeta}) = \widetilde {{\rho}_{j}}
(z^{\prime},g(z^{\prime},z^{\prime   \prime    \prime}),z^{\prime
\prime \prime},{\zeta})$.  For $  j  =
1,...,m - d, \nu = 1,...,p - d $  we define the functions

 $$\varphi^j_{\nu}(z^{\prime},z^{\prime \prime \prime},
\zeta) = {\partial  \widetilde  {\rho}^{j}  \over  \partial
z_{\nu}} - {\sum ^{d}_{t = 1} \widetilde a_{t {\nu}}
{\partial \widetilde {\rho}^{j} \over \partial z_{n -  d  +
t}}}, \tag{3.9} $$
where

$${\widetilde a}_{t{\nu}} = \sum^{d}_{s=1} {\widetilde b}_{ts}
{{\partial {\widetilde {\rho}^{s +  m  -  d}}}  \over  {\partial
z_{\nu}}},$$
and  $ (\widetilde b_{ts}  )$  is the inverse of the matrix

 $$  \left ( {{\partial {\widetilde {\rho}^{k}}} \over  {\partial
z_{l}}} \right )  ^{k = m - d + 1,...,m}_{l = N - d + 1,...,N}$$

Then the functions \thetag{3.9} are holomorphic on $ W(E,D) \times U$  and continuous
on $  ( W(E,D) \cup  E  )  \times  U$,   where $ U$  is a
neighborhood
of the origin in ${\Bbb  C}^{N}$ .
\par     Let $  (j(n),{\nu}_{1}(n),...,{\nu}_{k(n)}(n)), n = 1,...,s$
be the collections  of  positive  integers  from  Lemma 2.4.  We
define a set  $X \subset {\Bbb  C}^{2N}$     of the form

$$ X = X_{0}  \cap  X_{1}  \cap  ...  \cap  X_{s}  \cap  M^{c}  ,
\tag{3.10} $$
where

$$  X_{0}  = \lbrace (z,{\zeta}) \in {\bold C}^{2N}   :
 z^{\prime \prime} = g(z^{\prime},z^{\prime \prime \prime})
\rbrace, $$
and for  $n = 1,...,s$

 $$  X_{n}  = \lbrace (z,{\zeta}) \in {\bold C}^{2N}   :
{\varphi}^{j(n)}_{\nu}(z^{\prime},z^{\prime \prime \prime},
{\zeta}) = 0,
 {\nu}  = {\nu}_{1}(n),...,{\nu}_{k(n)}(n) \rbrace. \tag{3.11}$$

\proclaim { Lemma 3.2} X    can be represented in  the
following  form ( near the origin )
$$
\aligned
& z^{\prime \prime} = g(z^{\prime},z^{\prime \prime \prime}),\\
&{\zeta}_{\nu}    =   {\psi}_{\nu}(z^{\prime},z^{\prime   \prime
\prime},{\zeta}^{\prime}),  {\nu} = p - d + 1,...,N, \\
\endaligned
\tag{3.12}
$$
Here   $  {\zeta}^{\prime}  =  ({\zeta}_{1},...,{\zeta}_{p - d})$
and  the  functions   ${\psi}_{\nu}$    are holomorphic on
$W(E,D) \times U^{\prime}$   and of class $  C^1$    on $ (W(E,D)
\cup E) \times  U^{\prime}$,  where  $ U^{\prime}$   is a neighborhood
of the origin in ${\Bbb  C}^{p - d}$ .
\endproclaim
\par   {\it  Proof}. Let us consider the set $ X_{1}  \cap  ...  \cap
X_{s} \cap M^{c}$    that  is
defined by  equations  \thetag{3.11}  (  with  $n  =  1,...,s$  )  and
\thetag{3.8}.   We
compute the  Jacobian matrix $ J$  of  this  united system  ( which
contains $ ( N - p + d )$ equations ) with respect to the variables
${\zeta}_{p - d + 1},...,{\zeta}_{N}$   at  the origin. One has

$$ {\partial {\widetilde {{\rho}^{j}}} \over \partial z_{\nu}}
(z^{\prime},z^{\prime \prime \prime},{\zeta}) =
( {\partial {\widetilde {{\rho}_{j}}} \over
\partial z_{\nu}}) (z,{\zeta}) + \sum^{N-d}_{q = p - d  +
1} ( {\partial {\widetilde {{\rho}_{j}}}  \over  \partial
z_{q}}) (z,{\zeta})  {\partial  g_{q}  \over  \partial
z_{\nu}},$$
where  we   set    $z   =   (   z^{\prime},g(z^{\prime},z^{\prime
\prime \prime}),z^{\prime \prime \prime})$. We recall that  $g(0)
= 0$   and   ${\partial g_{q}} / {\partial z_{\nu}(0) = 0}$   for
any  $ q,{\nu}$.  Note also that $  p - d <
N - m + 1 $  ( see section 2).  Hence, \thetag{2.1}  implies
${\partial {\widetilde {\rho}_{j}}} / {\partial  z_{\nu}}(0)   =
0$  for  $ {\nu} = 1,...,p - d $  and each $ j$.  Moreover,
${\partial {\widetilde {{\rho}_{j}}}} / \partial z_{N - d + t}
(0) = 0,  j = 1,...,m - d,  t = 1,...,d$.  Therefore,

$${{\partial {\varphi}^{j}_{\nu}} \over {\partial {\zeta}_{\mu}}}
(0) = {{\partial^{2} {\widetilde  {{\rho}_{j}}}}  \over  {\partial
z_{\nu} \partial {\zeta}_{\mu}}}(0) =
{{{\partial^{2} {\rho}_{j}}} \over {\partial z_{\nu} \partial
{\overline z_{\mu}}}}(0) $$
for $ {\nu} = 1,...,p - d,{\mu} = p - d + 1,...,N - m, j =
1,...,  m  -
d$. Hence, $  J $  has the form

$$  J = \pmatrix H && \cr 0 & I_m \endpmatrix,$$
where $  H$   is the matrix (2.4), $ I_{m}$   is the identity
$ (m \times m)$-matrix.
     Now \thetag{2.3} implies $ det J \neq 0$.  By  the  implicit  function
theorem we get our statement, {\it Q.E.D}.
\par We  denote  by  $  X_E  $  the  set  defined  by  the   equations
\thetag{3.12}
under    the    condition    $    (z^{\prime},z^{\prime    \prime
\prime},{\zeta}^{\prime}) \in E \times U^{\prime}$.
\proclaim {Lemma 3.3} One has  the inclusion
$\widehat A_E  \subset  X_E $.
\endproclaim
\par   {\it  Proof}. Since  $ A_E  = \lbrace z \in {\Bbb  C}^{N} :
z^{\prime \prime} = g(z^{\prime},z^{\prime \prime \prime}),
(z^{\prime},z^{\prime \prime \prime}) \in E \rbrace $,  one can conclude that
$ (\partial {\widetilde \rho}^{j} / {\partial z_{\mu}}) \vert
\lbrace {\zeta} = {\overline z} \rbrace  = {\partial {\rho}^{j}} /
{\partial z_{\mu}}$  for  $ z  \in  A_E$ .
Hence, \thetag{3.9} and Lemma 3.1 imply
${\varphi}^{j}_{\nu} \vert \lbrace  {\zeta} = \overline z \rbrace =
0$   for $ z \in A_E$
. Since  $A_E  \subset M$, one has  ${\widetilde \rho}^{j} \vert
{\widehat A_E} = 0$.
Taking  into  account  the
equivalence of \thetag{3.10} and \thetag{3.12}, we obtain our  assertion,  {\it
Q.E.D}.

\head
4. THE REFLECTION PRINCIPLE.
\endhead

     We  consider  the  coordinates   $  (z,{\zeta})  \in  {\Bbb
C}^{2N} = {\Bbb  C}^{N} \times {\Bbb  C}^{N}$ .   Set   $
z = (\alpha,\beta)$,  where ${\alpha} = (z_{1} ,...,z_{n} )$,
${\beta} = (z_{n + 1},...,z_{N})$.
\proclaim {Lemma 4.1} Let  $ S$   be a generic $ C^1$ -manifold
in  ${\Bbb  C}^{n}_{\alpha}$   and
let $  W $  be  a  wedge  in ${\Bbb  C}^{n}$    with   the  edge $S$.
Suppose that $ F: W \to {\Bbb  C}^{N - n}_{\beta}  \times  {\Bbb
C}^{N}_{\zeta}$    is a holomorphic mapping of class  $C^1$
on  $ W \cup S$.  Let  $ Y \subset {\Bbb  C}^{2N}$     be the  graph  of
$ F$.  Assume  that
there exists a real  $n$-dimensional $ C^1$ -manifold
$ S^{\prime} \subset {\Bbb  C}^{2N}   = {\Bbb  C}^{N}_{z}  \times
{\Bbb  C}^{N}_{\zeta}$  such that  $ S^{\prime} \subset Y_S  \cap
\lbrace {\zeta} = {\overline z} \rbrace $ ( here
$ Y_S  = \lbrace (z,{\zeta}) : {(\beta},{\zeta}) = F({\alpha}),
{\alpha} \in S \rbrace )$. Then $  F $  extends holomorphically to
a neighborhood of $ S$.
\endproclaim
\par   {\it  Proof}. One  can  assume   $ 0 \in S^{\prime}$.  We  make
a  change  of coordinates of the form $ z = {\xi}^{\prime} + i
{\xi}^{\prime  \prime}, {\zeta}  =  {\xi}^{\prime}  -  i{\xi}^{\prime
\prime}$    ( here $ {\xi}^{\prime} = ({\xi}_{1},...,{\xi}_{N}),
{\xi}^{\prime \prime} =({\xi}_{N + 1},...,{\xi}_{2N}))$.
In the new coordinates the
diagonal $ \lbrace {\zeta} = {\overline z} \rbrace$  coincides with
${\Bbb  R}^{2N} = \lbrace {\xi} = ({\xi}^{\prime},  {\xi}^{\prime
\prime}): {\Im} {\xi}  = 0 \rbrace$.   Since   $  T_{0}{\overline
Y}$   is a complex n-dimensional  linear  space
in ${\Bbb  C}^{2N}$  ,  there exists a  n-dimensional coordinate
plane  ${\Pi}$  in  ${\Bbb  C}^{2N}$   such that
${\pi} : T_{0}{\overline Y} \to {\Pi}$   is an isomorphism  (
here  ${\pi} : {\Bbb  C}^{2N} \to {\Pi}$  is the natural projection ). One
can  assume   that  ${\Pi}$   is  the  plane  of  the   variables
${\xi}_{1},...,{\xi}_{n}$  . The  restriction ${\pi} : Y
\to {\Pi}$  is a local biholomorphism and ${\widetilde W} =
{\pi}(Y)$   is a wedge  in  ${\Pi}$    with  the   edge   $  Q  =
{\pi}(Y_S)$  that is generic in a  neighborhood of the origin.
Since $ S^{\prime} \subset {\Bbb  R}^{2N}$    and the restriction
${\pi} \vert S^{\prime}$
is a local diffeomorphism,  ${\pi}(S^{\prime})$   coincides with the plane  of
the variables  ${\Re} {\xi}_{1},..., {\Re} {\xi}_{n}$     in  a
neighborhood  of  the origin.  We denote this plane   by  $ {\Bbb
R}^n$ .   Since  $  S^{\prime}  \subset  Y_S $ , we
conclude that ${\Bbb  R}^{n}  \subset Q$.
\par     Let  $ G$   be a mapping holomorphic on $\widetilde W $
and  of  class  $C^1$  on ${\widetilde W} \cup Q$  so that $ Y$   is
the graph of  $ G$.  Then  $ S^{\prime}$ is
the graph of the restriction $  G \vert {\Bbb  R}^{n}$ .  Since
$ S^{\prime}  \subset  {\Bbb  R}^{2N}$  ,   the
restriction  $ G \vert {\Bbb  R}^{n}$    is a real valued mapping.
Therefore,  one
can apply the edge of the wedge theorem \cite{Ru} to the mappings $ G$
and $ G^{\ast}  = \overline{ G}(\overline{\xi}_{1},...,
\overline{\xi}_{n})$. Thus, $ G$  extends  holomorphically  to  a
neighborhood of the origin. Therefore, $ Y$  extends to  a  complex
manifold in a  neighborhood  of  the  origin.  This  implies  our
assertion, {\it Q.E.D}.
\par  {\it Completion of the proof of Theorem 1}. Setting  $ S = E
\times U^{\prime}$  and  $ S^{\prime} = {\widehat A_E}$ ,  we apply
Lemma  4.1   to   the   mapping   $  F  =  (g,{\psi}_{p  -   d   +
1},...,{\psi}_{N})$  ( Lemmas 3.2 and 3.3 show that one can do it
).
We get that $ F$   extends holomorphically to  a  neighborhood  of
the origin. Hence,$ g$  extends  holomorphically.  But  this  means
that $ f$  extends holomorphically to a neighborhood of the  origin
in ${\Bbb  C}^{p}$ . Thus, the mapping  ${\Phi}(x) = (x,f(x))$   is
holomorphic  near
the origin and the condition (ii) of Definition  1  implies  that
$E = \lbrace x: {\rho}_{j} \circ {\Phi}(x) = 0, j = m - d + 1,...
,m \rbrace $. Hence,$  E$  is real analytic,{\it Q.E.D}.

\head
5. CERTAIN SPECIAL CASES.
\endhead

   {\it  Proof of Proposition 2}. Without loss of generality we assume
$a = 0,  a^{\prime} = f(a) = 0$. It is sufficient to show  that
$ A =  {\Gamma}_{f}$  and $  M = S \times S^{\prime}$   are Levi
transverse at the origin. Let $ r^{\prime}_{j} , j
= 1,...,d'$  be the defining functions of $ S^{\prime}$  and  let
$  r_{j} ,  j  = 1,...,d$  be the defining functions of $ S$.
Considering  $ {\rho}_{j}  = r^{\prime}_{j}  -
\sum {\lambda}_{j{\nu}}r_{\nu}$ ,    we  get  the  defining  functions  of
$M$    with  the condition   $d{\rho}_{j} \vert T_{0} A = 0$. The
matrices of the  Levi  operators $ L^{j}_{0}$
of the functions  ${\rho}_{j}$  ( on $ T^{c}_{0}M$ ) have the form

$$ \pmatrix R^{j} & 0 \cr 0 & L^{{\prime}^j}_{0} \endpmatrix$$
where  $ R^j$  being the  operators on $ T^{c}_{0}S$, and
$L^{{\prime}^{j}}_{0}$ being the Levi operators of $r^{\prime}_{j}$
( we identify
operators  and their matrices ).  Since  $  E  =  S  $   and    $
T^{c}_{0}A_E  = \lbrace (v,df_{0}(v)):  v \in T^{c}_{0}S \rbrace$,
one can conclude by \thetag{1.13} that the  dimension  of  the
space $ L = \sum^{d^{\prime}}_{j = 1} L^{j}_{0}(T^{c}_{0}A_E) $ is equal
to $ p^{\prime} -  d^{\prime}$.  Since $  L $  and $
T^{c}_{0}(A_E)$  are orthogonal by Lemma 2.3,  we  get  \thetag{1.9},{\it
Q.E.D}.
\par     Let $ < , >$  be   a   hermitian   scalar   product   on
${\Bbb C}^{n^{\prime}}$     defining
the Levi operators. Assume that the  conditions  of  Corollary  1
hold. We show that in this case \thetag{1.13} also  holds.  Assume  that
this is not true. Then there exists a vector ${\xi} \in T^{c}_{0}
S^{\prime} \backslash \lbrace 0 \rbrace $
orthogonal  (  with  respect  to  $  < ,  >$   )   to   the   space
$\sum^{d^{\prime}}_{j                    =                     1}
L^{{\prime}^j}_{0}(df_{0}(T^{c}_{0}S))$.   Since  $ df_{0}$    is
surjective,        one       has       $df_{0}(T^{c}_{0}S)      =
T^{c}_{0}S^{\prime}$.                 Therefore,                $
<L^{{\prime}^j}_{o}({\xi}),{\eta}>                              =
\overline {<L^{{\prime}^j}_{0}({\eta}),{\xi} >} = 0$  for
each  $ {\eta} \in T^{c}_{0}S^{\prime},  j = 1,...,d^{\prime}$
( $L^{{\prime}^j}_{0}$   is hermitian ). We obtain
a  contradiction  to the  condition  of  Levi  non-degeneracy  of
$ S^{\prime}$.  This   proves Corollary 1.
\par     {\it Proof  of  Corollary  4}. The condition (i) of  the Levi
transversality
is trivial. Since the intersection  $  T_{0}  A  \cap
T_{0}M$   is a totally  real  space in $ T_{0}A$, the  condition  (ii)
holds.   Since $  T^{c}_{0}A_E   =   T^{c}_{0}M   =   \lbrace    0
\rbrace,$   the   condition   (iii)   also   is
trivial.  In  this  case  we  do  not  need       the   functions
\thetag{3.9}   to   define   the   set  $  X$.   Hence,   the   proof   of
Theorem 1 is valid for $C^1$ - wedges and $C^1$ -mappings in this special
case, {\it Q.E.D}.
\par     Corollary 2 is a special case of Corollary 4; the  proof  of
Corollary 3 is quite similar.

\head
6. PROOF OF THEOREM 2.
\endhead

     Set $ E = \lbrace x \in U : r_{j}(x) = 0, j = 1,...,n
\rbrace $, $ \partial r_{1} \land ... \land \partial r_{n} \neq
0$, $  E^{\prime} = \lbrace x^{\prime} \in V :
{\psi}_{j}(x^{\prime},\overline   x^{\prime})   =    0,    j    =
1,...,n^{\prime} \rbrace  $,  $\partial  {\psi}_{1}  \land   ...   \land
\partial {\psi}_{n^{\prime}} \neq 0 $.
 Let  $\varphi(x,\overline x)$  be a real analytic defining  function  of
$b{\Omega}$.  Without loss of generality we assume that $  0  \in
E \cap U,  f(0) = 0,   r_{j}   =  x_{j}   +  \overline  x_{j}   +
o(\vert x \vert),    {\psi}_{j}   =  x^{\prime}_{j}  +  \overline
x^{\prime}_{j} + o(\vert x^{\prime} \vert),   {\varphi} = x_{n}  +
\overline x_{n}  + \sum^{n - 1}_{j = 1} \vert x_{j} \vert^{2}  +
o({\vert x \vert}^2 )$.
\par     Since  $ f(E)  \subset  E^{\prime}$    and  $  {\rho}  \vert
E^{\prime} = 0$,   we have $ ({\rho} \circ f)(x) =
 {\rho}(f(x),\overline f(x)) = 0$   for $ x \in E$.  By the division
lemma   one  can conclude that

$$ {\rho} \circ f = \sum^{n-1}_{j = 1} {\lambda}_{j}r_{j}(x) +
{\lambda}_{n} {\varphi}(x,\overline x), \tag{6.1}$$
for $ x \in U$.  Here  $ {\lambda}_{j} \in C^1(U)$.
\proclaim  { Lemma 6.1}      We have

$$ {\lambda}_{j}(x) = 0,  j = 1,...,n - 1,  {\lambda}_{n}(x) \neq
0$$
for $  x \in E \cap U$.
\endproclaim
\par     {\it   Proof}.   Since  $  f({\Omega}   \cap    U)    \subset
{\Omega}^{\prime} \cap V$,  we get that $ {\rho} \circ f $   is  a
negative plurisubharmonic function on ${\Omega} \cap U$.  By   Hopf's  lemma
we  conclude   $  \vert   {\rho}   \circ   f(x)   \vert   \geq   C
dist(x,b{\Omega})$   for each $ x \in {\Omega} \cap U $  and
a positive constant  C  ( here $dist$ is the Euclidean distance ).
Hence,$ d({\rho} \circ f)(0) \neq 0$. Since the domain
${\Omega} \cap U = \lbrace x \in U : \varphi(x,\overline x)
< 0 \rbrace $  is contained in the domain $  D = \lbrace x \in U :
({\rho} \circ f) < 0  \rbrace$,
the tangent planes of their boundaries coincide at each point of
$E  \subset  b{\Omega}  \cap  bD$.  This   implies   the   desired
statement,{\it Q.E.D}.
     Lemma 6.1 and \thetag{6.1} imply

$$ {{\partial{}} \over {\partial x_{\nu}}}({\rho} \circ f)(x) =
{\lambda}_{n}(x) {{\partial{}} \over \partial x_{\nu}}{\varphi}(x)
, {\nu} = 1,...,n , \tag{6.2}$$
for $ x \in E \cap U$. We get

$${\partial{\varphi} \over {\partial x_{n}}}
{{\partial({\rho} \circ f)} \over {\partial x_{\nu}}}
 - {\partial {\varphi} \over {\partial x_{\nu}}} {{\partial({\rho}
\circ f)} \over {\partial x_{n}}} = 0,  {\nu} = 1,...,n - 1,
\tag{6.3}$$
for  $ x \in E \cap U$.
     Set    $y  \in  {\Bbb   C}^{n}$  ,   $y^{\prime}  \in  {\Bbb
C}^{n^{\prime}}  ,  {\Bbb  C}^N   =  {\Bbb  C}^n   \times  {\Bbb
C}^{n^{\prime}}$ .  We  define  the functions

 $$ h_{\nu}(x,y,y^{\prime}) = {\partial {\varphi} \over {\partial
x_{n}}}(x,y){\partial{\rho}(f(x),y^{\prime})   \over    {\partial
x_{\nu}}} - {\partial {\varphi} \over {\partial x_{\nu}}}(x,y)
{\partial{\rho}(f(x),y^{\prime}) \over {\partial x_{n}}},
{\nu} = 1,...,n - 1, \tag{6.4}$$

     We consider  the  set   $A$   in $ {\Bbb   C}^{2N}  =  {\Bbb
C}^{n}_{x}  \times  {\Bbb    C}^{n^{\prime}}_{x^{\prime}}   \times
{\Bbb  C}^{n}_{y} \times {\Bbb  C}^{n^{\prime}}_{y^{\prime}}$     that
is defined by the equations

$$     h_{\nu}(x,y,y^{\prime}) = 0, {\nu} = 1,...,n - 1,$$
$$  \varphi(x,y) = 0,  \tag{6.5}$$
$$ {\psi}_{j}(x^{\prime},y^{\prime}) = 0, j = 1,...,n^{\prime},$$
and

$$ x^{\prime} = f(x). \tag{6.6}$$

     Here  $  x  \in  {\overline  {\Omega}}  \cap  U  $   and   $
(y,y^{\prime}) $  belongs to a  neighborhood
of the origin in $ {\bold C}^N$ .
\proclaim {Lemma 6.2} The  set  $ A$   defined  by \thetag{
6.5}, \thetag{6.6},   is  a
complex   n-dimensional  manifold   with    $ C^1$ -boundary   in   a
neighborhood of the origin in  ${\Bbb C}^{2N}$ .
\endproclaim
\par {\it Proof}. We compute the Jacobian matrix $ J $ of \thetag{6.5}  with
respect to the variables $ (y,y^{\prime}) $ at $ 0$.  For
$ {\mu} =  1,...,n  -  1$  we have

$${\partial h_{\nu} \over {\partial y_{\mu}}}(0) =
{\partial^{2}{\varphi} \over {\partial x_{n} \partial y_{\mu}}}
(0,0)   {\partial {\rho}(f(x),y^{\prime})    \over    {\partial
x_{\nu}}}(0,0) +
{\partial           {\varphi}           \over           {\partial
x_{n}}}(0,0)
{\partial^{2}{\rho}(f(x),y^{\prime})                         \over
{\partial{y_{\mu}} \partial{x_{\nu}}}}(0,0) -$$
$$-{\partial^{2}{\varphi}   \over   {\partial    x_{\nu}    \partial
y_{\mu}}}(0,0) {\partial {\rho}(f(x),y^{\prime}) \over  {\partial
x_{n}}}(0,0) -
{\partial{\varphi} \over {\partial x_{\nu}}}(0,0)
 {\partial^{2}{\rho}(f(x),y^{\prime})  \over   {\partial   y_{\mu}
\partial {x_n}}}(0,0).$$

     One has

$${\partial^{2}{\varphi}   \over    {\partial    x_{n}    \partial
y_{\mu}}}(0,0) =
{\partial^{2}{\rho}(f(x),y^{\prime})  \over   {\partial   y_{\mu}
\partial  x_k}}(0,0)   =   {\partial{\varphi}   \over   {\partial
x_{\nu}}}(0,0) = 0,$$
for  $ {\mu},{\nu} = 1,...,n - 1, k = 1,...,n$.
     Also,
$${\partial^{2}{\varphi}   \over    {\partial    x_{\nu}
\partial y_{\mu}}}(0,0) = {\delta}_{{\nu}{\mu}}$$
    ( the Kronecker symbol  ).  Lemma 6.1 implies
$${\partial {\rho}(f(x),y^{\prime})
\over {\partial x_{n}}}(0,0) = {\alpha} \neq 0$$.
  Thus,

$$ J = \pmatrix -{\alpha}I_{n - 1} & & \cr 0 & & 1 & \cr
0 & & 0 & & I_{n^{\prime}} \endpmatrix$$
and $  det J \neq 0 $. Hence, the implicit  function  theorem
implies
that  $A$  can be represented in the following form :

$$x^{\prime} = f(x),  y = g(x),  y^{\prime} = p(x). \tag{6.7}$$
and the mapping  $(f,g,p)$  is holomorphic on ${\Omega} \cap  U $
and  of
class  $ C^1$    on  $\overline{\Omega} \cap U$,{\it  Q.E.D}.
\par     Let $W$   be a wedge in  ${\Omega \cap  U}$   with  the  edge
$E$. In view of \thetag{6.7} one can assume that $ A$    is  the  graph
of  the  mapping  $F = (f,g,p)$  over $  W$.
\proclaim { Lemma 6.3} The set
$$     A_E       =      \lbrace      (x,y,x^{\prime},y^{\prime}):
(x^{\prime},y,y^{\prime}) = F(x), x \in E \rbrace $$
is contained in the diagonal

$$  M = \lbrace (x,x^{\prime},y,y^{\prime}) : x  =  \overline  y,
x^{\prime}  = \overline y^{\prime} \rbrace.$$
\endproclaim
\par     {\it Proof}.  We  set  $   x  =\overline  y,  x^{\prime}   =
\overline y^{\prime},  x  \in  E $  in  \thetag{6.5}, \thetag{6.6}.
 Since  $  E  \subset  b{\Omega}$    and  $   f(E)  \subset
E^{\prime}$,  we get $  {\varphi}(x,\overline  x)   =  0  $   and
${\psi}_{j}(x^{\prime},\overline x^{\prime})   =   0$.   In   view  of
\thetag{6.3} and  \thetag{6.4}
  we  obtain  $h_{\nu}(x,\overline x, \overline x^{\prime})  =
h_{\nu}(x, \overline x, \overline f(x))  =  0$.  Hence,  the  equivalence   of
\thetag{6.5}, \thetag{6.6} and \thetag{6.7} implies that $A_E$ and $A_E  \cap  M$  have
the same real dimension. This gives the desired assertion, {\it
Q.E.D}.
\par     Thus, one can apply Corollary 4 to  $A$  and $ M$.  We conclude
that the mapping  $F = (f,g,p)$   ( and, certainly,$  f $)  extends
holomorphically to a neighborhood of the origin and $  E $  is real
analytic, {\it Q.E.D}.

\head
Alexander Sukhov
Department of Mathematics
Bashkir State University
450074, Frunze str. 32,
Ufa, Russia
\endhead
\end
\Refs
\widestnumber\key{AAAA}
\ref\key{Al} \by H.Alexander \paper Continuing   1-dimensional
analytic sets \jour Math.Ann \vol  191 \yr 1971 \pages 143-144
\endref
\ref\key{BJT} \by  M.Baouendi, H.Jacobowitz, F.Treves \paper On the
analyticity of CR mappings \jour Ann. Math \vol 122 \yr 1985
\pages 365-400 \endref
\ref\key{Ch1} \by E.M.Chirka \paper Regularity of boundaries  of  analytic
sets \jour Matem. sb \vol 117 \yr 1982 \pages 291-336 \endref
\ref\key{Ch2} \by  E.M.Chirka \paper Introduction to the geometry of
CR manifolds \jour Uspekhi matem.nauk \vol 46 \yr 1991 \pages 81-164
\endref
\ref\key{Der} \by  M.Derridj Le principe de  reflexion  en
des  points  de faible  pseudoconvexite,  pour  des
applications   holomorphes propres \jour Invent. math \vol 79
\yr 1985 \pages 197-215 \endref
\ref\key{DW} \by  K.Diederich ,  S.Webster \paper  A  reflection
principle  for degenerate real hypersurfaces \jour  Duke Math.J
\vol 47 \yr 1980 \pages 835-843 \endref
\ref\key{DF} \by    K.Diederich , J.E.Fornaess \paper  Proper holomorphic
mappings between real analytic pseudoconvex domains in ${\Bbb
C}^{n}$ \jour  Math.Ann \vol 282 \yr 1988 \pages 681-700 \endref
\ref\key{Fe} \by C.Fefferman \paper The Bergman kernel and  biholomorphic
mappings of pseudoconvex domains \jour Invent. Math \vol 26 \yr 1974
\pages 1-65 \endref
\ref\key{Fo} \by F.Forsneric \paper  On the boundary regularity  of  proper
mappings \jour Ann. Scuol. Norm \vol 13 \yr 1986 \pages 109 - 128
\endref
\ref\key{Le} \by  H.Lewy \paper On  the  boundary  behaviour   of
holomorphic mappings \jour  Acad. Naz. Lincei \vol  35  \yr  1977
\pages 1-8 \endref
\ref\key{NR} \by A.Nagel, J.-P.Rosay \paper Maximum   modulus   sets   and
reflection sets \jour Ann. Inst. Fourier \vol 41 \yr 1991 \pages 431-466
\endref
\ref\key{Pi} \by S.Pinchuk \paper On the analyic continuation  of
biholomorphic mappings \jour Matem. Sb \vol  98 \yr 1975
\pages 416-435 \endref
\ref\key{PCh} \by  S.Pinchuk, E.Chirka \paper On  the  reflection  principle
for analytic  sets \jour  Izv.  Acad.  Nauk  USSR,Ser.Matem \vol 52
\yr 1988 \pages 200-210 \endref
\ref\key{Sh} \by  B.Shiffman \paper  On the continuation of analytic
curves \jour Math.Ann \vol 184 \yr 1970 \pages 268-274 \endref
\ref\key{Su1} \by   A.Sukhov \paper On the analytic  continuation
across  generic manifolds \jour  Matem. zametki \vol 46 \yr 1989
\pages 119-120 \endref
\ref\key{Su2} \by  A.Sukhov \paper On the extension  of  analytic
sets  through generic manifolds \jour Russian Acad. Sci. Izv. Math.
\vol 401 \yr 1993 \pages 203 - 212 \endref
\ref\key{Su3} \by  A.Sukhov \paper On   algebraicity   of  complex   analytic
sets \jour Math. USSR Sbornik \vol 74 \yr 1993 \pages 419 - 426
\endref
\ref\key{Su4} \by  A.Sukhov \paper On holomorphic mappings of domains of the
type "wedge" \jour Matem. zametki \vol 52 \yr 1992 \pages 141-145
\endref
\ref\key{Su5} \by  A.Sukhov \book The reflecion principle and a
mapping  problem for Cauchy-Riemann manifolds \publ preprint
\yr 1992  \publaddr Univ.  Provence,Marseille \endref
\ref\key{Su6} \by A.Sukhov \paper  On CR mappings of real
quadric  manifolds \jour
Mich. Math. J \vol 41 \yr 1994 \pages 143 -150 \endref
\ref\key{Ru} \by   W.Rudin \book Lectures on   the   Edge   of   the   Wedge
theorem \publ
Regional Conf.Series in  Math.  A.M.S.  \publaddr  Providence,  R.I.,6
(1971) \endref
\ref\key{TH} \by  A.Tumanov  ,   G.Henkin \paper  Local
characterization   of
holomorphic automorphisms  of  Siegel  domains  \jour  Func.
Anal and Appl \vol 17 \yr 1983 \pages 49-61 \endref
\ref\key{We1} \by S.Webster \paper On the reflection principle in
several  complex
variables \jour  Proc. Amer. Math. Soc \vol 71 \yr 1978 \pages 26-28
\endref
\ref\key{We2} \by S.Webster \paper Holomorphic mappings  of  domains
with  generic corners \jour Proc. Amer. Math. Soc \vol 86  \yr 1982
\pages 236-240 \endref
\endRefs
